\documentclass[12pt,twoside]{article}

\usepackage[english]{babel}
\usepackage{umj_1}           

\usepackage{graphicx, srcltx}   
\usepackage{cite}               

\usepackage{amsmath}            
\usepackage{amsfonts,amssymb}   
\usepackage{srcltx}             

\newcommand{\Leftclt}{Roman A. Veprintsev}
\newcommand{\Rightclt}{Dunkl harmonic analysis and fundamental sets of continuous functions}
\usepackage{authblk}

\begin{document}

\begin{Titul}
{\large \bf DUNKL HARMONIC ANALYSIS AND\\ FUNDAMENTAL SETS OF CONTINUOUS FUNCTIONS\\[0.2em] ON THE UNIT SPHERE }\\[3ex]
{{\bf Roman~A.~Veprintsev} \\[5ex]}
\end{Titul}

\begin{Anot}
{\bf Abstract.} We establish a necessary and sufficient condition on a continuous function on $[-1,1]$ under which the family of functions on the unit sphere $\mathbb{S}^{d-1}$ constructed in the described manner is fundamental in $C(\mathbb{S}^{d-1})$. In our construction of functions and proof of the result, we essentially use Dunkl harmonic analysis.


\end{Anot}


\section{Introduction and preliminaries}

We need some elements of the general Dunkl theory (see \cite{dunkl_xu_book_orthogonal_polynomials:2014,dunkl_article_integral_kernels:1991,dunkl_article_reflection_groups:1988,dai_xu_book_approximation:2013,dunkl_article_operators:1989}); for a background on reflection groups and root systems the reader is referred to \cite{Humphreys_book_reflection:1990,dunkl_xu_book_orthogonal_polynomials:2014}.

Let $\mathbb{R}^d$ denote $d$-dimensional Euclidean space. For $x\in\mathbb{R}^d$, we write $x=(x_1,\dots,x_d)$. The inner product of $x,y\in\mathbb{R}^d$ is denoted by $\langle x,y\rangle=\sum\limits_{i=1}^d x_iy_i$, and the norm of $x$ is denoted by $\|x\|=\sqrt{\langle x,x\rangle}$.

The unit sphere $\mathbb{S}^{d-1}$, $d\geq 2$, and the unit ball $\mathbb{B}^d$ of $\mathbb{R}^d$ are defined by
\begin{equation*}
\mathbb{S}^{d-1}=\{x\colon\, \|x\|=1\}\quad\text{and}\quad \mathbb{B}^d=\{x\colon\, \|x\|\leq 1\}.
\end{equation*}

For a nonzero vector $v\in\mathbb{R}^d$, define the reflection $\sigma_v$ by
\begin{equation*}
\sigma_v(x)=x-2\frac{\langle x,v\rangle}{\|v\|^2}\,v,\quad x\in\mathbb{R}^d.
\end{equation*}
Each reflection $\sigma_v$ is contained in the orthogonal group $O(\mathbb{R}^d)$.

We give some basic definitions and notions which will be important.

\begin{definition}
Let $R\subset\mathbb{R}^d\setminus\{0\}$ be a finite set. Then $R$ is called a root system if

$(1)$ $R\cap \mathbb{R}v=\{\pm v\}$ for all $v\in R$;

$(2)$ $\sigma_v(R)=R$ for all $v\in R$.

The subgroup $G=G(R)\subset O(\mathbb{R}^d)$ which is generated by the reflections $\{\sigma_v\colon\, v\in R\}$ is called the reflection group associated with $R$.
\end{definition}

For any root system $R$ in $\mathbb{R}^d$, the reflection group $G=G(R)$ is finite. The set of reflections contained in $G(R)$ is exactly $\{\sigma_v\colon\, v\in R\}$.

Each root system can be written as a disjoint union $R=R_+\cup -R_+$, where $R_+$ and $-R_+$ are separated by a hyperplane through the origin. Such a set $R_+$ is called a positive subsystem. Its choice is not unique.

\begin{definition}
A nonnegative function $\kappa$ on a root system $R$ is called a multiplicity function on $R$ if it is $G$-invariant, i.e. $\kappa(v)=\kappa(g(v))$ for all $v\in R$, $g\in G$.
\end{definition}

\begin{definition}
The Dunkl operators are defined by
\begin{equation*}
\mathcal{D}_if(x)=\frac{\partial f(x)}{\partial x_i}+\sum\limits_{v\in R_+} \kappa(v)\frac{f(x)-f(\sigma_v(x))}{\langle x,v\rangle} \langle v,e_i\rangle,\quad 1\leq i\leq d,
\end{equation*}
where $e_1,\dots,e_d$ are the standard unit vectors of $\mathbb{R}^d$.
\end{definition}

The above definition does not depend on the special choice of $R_+$, thanks to the $G$-invariance of $\kappa$. In case $\kappa=0$, the Dunkl operators reduce to the corresponding partial derivatives.

Suppose $\Pi^d$ is the space of all polynomials in $d$ variables with complex coefficients, $\mathcal{P}_n^d$ $\bigl(n\in \mathbb{N}_0=\{0,1,2,\dots\}\bigr)$ is the subspace of homogeneous polynomials of degree $n$ in $d$ variables.

According to \cite{dunkl_article_integral_kernels:1991}, there exists a unique linear isomorphism $V_\kappa$ of $\Pi^d$ such that
\begin{equation*}
V_\kappa(\mathcal{P}_n^d)=\mathcal{P}_n^d,\,\,\, n\in\mathbb{N}_0,\quad V_\kappa1=1,\quad\text{and}\quad\mathcal{D}_iV_\kappa=V_\kappa\frac{\partial}{\partial x_i},\quad 1\leq i\leq d.
\end{equation*}
This operator is called the Dunkl intertwining operator. If $\kappa=0$, $V_\kappa$ becomes the identity operator.

M.~R\"{o}sler has proved in \cite{rosler_article_positivity:1999} that for each $x\in\mathbb{R}^d$ there exists a unique probability measure $\mu_x^\kappa$ on the Borel $\sigma$-algebra of $\mathbb{R}^d$ with support in $\{y\colon\,\|y\|\leq\|x\|\}$, such that for all polynomials $p$ on $\mathbb{R}^d$ we have
\begin{equation}\label{rosler_measure}
V_\kappa p(x)=\int\nolimits_{\mathbb{R}^d} p(y)\,d\mu_x^\kappa(y).
\end{equation}

For $X=[-1,1]$ or $X=\mathbb{S}^{d-1}$, we denote by $C(X)$ the space of continuous complex-valued functions on $X$.

\begin{definition}
We define the Dunkl truncated intertwining operator at an arbitrary point $\xi\in\mathbb{S}^{d-1}$
\begin{equation*}
V_\kappa(\xi)\colon C([-1,1])\to C(\mathbb{S}^{d-1})
\end{equation*}
by
\begin{equation*}
V_\kappa(\xi;g,x)=\int\nolimits_{\mathbb{B}^d} g(\langle \xi,\zeta\rangle)\, d\mu_x^\kappa(\zeta),\quad x\in\mathbb{S}^{d-1},\quad g\in C[-1,1],
\end{equation*}
where $\mu_x^\kappa$ is the measure given in \eqref{rosler_measure}.
\end{definition}

The operator $V_\kappa(\xi)$ is well defined. Indeed, there exists a sequence $\{p_n\}$ of polynomials in $\Pi^1$ such that
\begin{equation*}
\sup\limits_{t\in[-1,1]} \, |g(t)-p_n(t)|\to 0,\quad n\to\infty.
\end{equation*}
Note that $p_n(\langle \xi,\cdot\rangle)\in\Pi^d$. Thus, for all $x\in\mathbb{B}^d$,
\begin{equation*}
\begin{split}
\bigl|V_\kappa&(\xi;g,x)-V_\kappa\bigl[p_n(\langle \xi,\cdot\rangle)\bigr](x)\bigr|\leq \int\nolimits_{\mathbb{B}^d} \bigl|g(\langle\xi,\zeta\rangle)-p_n(\langle\xi,\zeta\rangle)\bigr| \, d\mu_x^\kappa(\zeta)\\
&\leq \sup\limits_{\zeta\in\mathbb{B}^d} \, \bigl|g(\langle\xi,\zeta\rangle)-p_n(\langle\xi,\zeta\rangle)\bigr|=\sup\limits_{t\in[-1,1]} \, |g(t)-p_n(t)|\to 0,\quad n\to\infty.
\end{split}
\end{equation*}
As $V_\kappa\bigl[p_n(\langle\xi,\cdot\rangle)\bigr]$ is continuous on $\mathbb{R}^d$, then we deduce that $V_\kappa(\xi;g)$ is continuous on the unit ball $\mathbb{B}^d$.

In the present paper, we establish a necessary and sufficient condition on a function $g\in C([-1,1])$ under which the family of functions $\{V_\kappa(x;g)\colon x\in\mathbb{S}^{d-1}\}$ is fundamental in $C(\mathbb{S}^{d-1})$. This result generalizes Theorem 2 in \cite{sun_cheney_article_fundamental_sets:1997}. We state and prove the main result in section~\ref{section_with_main_result}.

Recall that a set $\mathcal{F}$ in a Banach space $\mathcal{E}$ is said to be fundamental if the linear span of $\mathcal{F}$ is dense in $\mathcal{E}$.

\section{Fundamentality in $C(\mathbb{S}^{d-1})$}

Suppose that the unit sphere $\mathbb{S}^{d-1}$ is equipped with a positive Borel measure and notions of orthogonality for functions on $\mathbb{S}^{d-1}$ are defined in terms of this Borel measure.

Let there be given an orthogonal sequence $\{U_n\}_{n=1}^\infty$ of finite-dimensional subspaces in $C(\mathbb{S}^{d-1})$. It is assumed that $\bigcup\limits_{n=1}^\infty U_n$ is fundamental in the space $C(\mathbb{S}^{d-1})$. For each $n$, let $\{u_{nj}\colon 1\leq j\leq\dim U_n\}$ denote a real-valued orthonormal basis of $U_n$.

Assume that we possess a summability method, given by an infinite matrix $A\!=\!\!\bigl[A_{nm}\bigr]_{n,m=1}^\infty$ with complex entries that has these properties:

\begin{itemize}
\item[(i)] each row of $A$ has only finitely many nonzero elements;

\item[(ii)] $\lim\limits_{n\to\infty} \, A_{nm}$ exists for each $m=1,2,\ldots$;

\item[(iii)] the sequence of functions $f_n(x,y)=\sum\nolimits_{m} \, A_{nm}\,\sum\nolimits_{j} \, u_{mj}(x)u_{mj}(y)$, where $x,y\in\mathbb{S}^{d-1},$ converges (as $n\to\infty$) uniformly in $x$ and $y$ to a limit function $f(x,y)$.
\end{itemize}

\begin{teoen}\label{main_general_theorem}
Suppose that the hypotheses given above are satisfied, and let $f$ be as in {\upshape (iii)}. In order that the set of functions $\{x\mapsto f(x,y)\colon\, y\in\mathbb{S}^{d-1}\}$ be fundamental in $C(\mathbb{S}^{d-1})$ it is necessary and sufficient that $\lim\limits_{n\to\infty} \, A_{nm}$ be nonzero for all $m$.
\end{teoen}

This theorem was proved in \cite[Section~2, Theorem~1]{sun_cheney_article_fundamental_sets:1997} in a more general setting.

\section{Some facts of Dunkl harmonic analysis on the unit sphere}\label{section_with_Dunkl_harmonic_analysis}

The Dunkl Laplacian is defined by
\begin{equation*}
\Delta_\kappa=\mathcal{D}_1^2+\dots+\mathcal{D}_d^2.
\end{equation*}
The Dunkl Laplacian plays the role of the ordinary Laplacian. In the special case $\kappa=0$, $\Delta_\kappa$ reduces to the ordinary Laplacian.

A $\kappa$-harmonic polynomial $P$ of degree $n\in\mathbb{N}_0$ is a homogeneous polynomial $P\in \mathcal{P}_n^d$ such that $\Delta_\kappa P=0$. The $\kappa$-spherical harmonics of degree $n$ are the restriction of $\kappa$-harmonics of degree $n$ to the unit sphere $\mathbb{S}^{d-1}$. Let $\mathcal{A}_n^d(\kappa)$ be the space of $\kappa$-spherical harmonics of degree $n$ and let $N(n,d)$ be the dimension of $\mathcal{A}_n^d(\kappa)$.

The weighted inner product of $f,h\in C(\mathbb{S}^{d-1})$ is denoted by
\begin{equation*}
\langle f,h\rangle_\kappa=\frac{1}{\sigma_d^\kappa}\int\nolimits_{\mathbb{S}^{d-1}} f(x)\overline{h(x)}\, w_\kappa(x)d\omega(x),
\end{equation*}
where $d\omega$ is the Lebesgue measure on $\mathbb{S}^{d-1}$, $w_\kappa$ is the weight function, invariant under the reflection group $G$, defined by
\begin{equation*}
w_\kappa(x)=\prod\limits_{v\in R_+} |\langle v,x\rangle|^{2\kappa(v)},\quad x\in\mathbb{S}^{d-1},
\end{equation*}
and $\sigma_d^\kappa$ is the constant chosen  such that $\langle 1,1\rangle_{\kappa}=1$.

Note that if $\kappa=0$, then $w_\kappa=1$.

In the rest of this section, we assume that $\kappa\not=0$ if $d=2$.

The following properties hold:

\begin{itemize}
\item[(I)] $\bigcup\limits_{n=0}^\infty \mathcal{A}_n^d(\kappa)$ is fundamental in $C(\mathbb{S}^{d-1})$.

\item[(II)] For any real-valued orthonormal basis $\{S_{n,j}^{\kappa}\colon\, 1\leq j\leq N(n,d)\}$ of $\mathcal{A}_n^d(\kappa)$,
    \begin{equation*}
    V_\kappa(x;C_n^{\lambda_\kappa},y)=\frac{\lambda_\kappa}{n+\lambda_\kappa}\, \sum\limits_{j=1}^{N(n,d)} S_{n,j}^{\kappa}(x)S_{n,j}^{\kappa}(y),\quad x,y\in\mathbb{S}^{d-1},
    \end{equation*}
    where $C_n^{\lambda_\kappa}(\cdot)$ denotes the Gegenbauer polynomial, $\lambda_\kappa$ is a positive constant defined by
    \begin{equation}\label{lambda_kappa}
    \lambda_\kappa=\sum\limits_{v\in R_+} \kappa(v)+\frac{d-2}{2}
    \end{equation}
(see, for example, \cite[Section~7.2]{dai_xu_book_approximation:2013}).
\item[(III)] If $n\not=m$, then $\mathcal{A}_{n}^{d}(\kappa)\perp\mathcal{A}_m^d(\kappa)$, i.e. $\langle P,Q\rangle_\kappa=0$ for $P\in\mathcal{A}_{n}^{d}(\kappa)$ and $Q\in \mathcal{A}_m^d(\kappa)$ \cite[Theorem~1.6]{dunkl_article_reflection_groups:1988}.
\end{itemize}

Property (I) follows from the Weierstrass approximation theorem: if $f$ is continuous on $\mathbb{S}^{d-1}$, then it can be uniformly approximated by polynomials restricted to $\mathbb{S}^{d-1}$. According to \cite[Theorem~1.7]{dunkl_article_reflection_groups:1988}, these restrictions belong to the linear span of $\bigcup\limits_{n=0}^\infty \mathcal{A}_n^d(\kappa)$.

\section{Main result and its proof}\label{section_with_main_result}

We can now establish the main result of the paper.

\begin{teoen}
Fix $d\geq 2$. Suppose $R$ is a fixed root system in $\mathbb{R}^d$, $\kappa$ is a multiplicity function on $R$. Assume that $\kappa\not=0$ if $d=2$. Let $g\in C([-1,1])$. In order that the family of functions $\{V_\kappa(x;g)\colon\, x\in\mathbb{S}^{d-1}\}$ be fundamental in $C(\mathbb{S}^{d-1})$ it is necessary and sufficient that
\begin{equation*}
\int\nolimits_{-1}^{1} g(t) C_n^{\lambda_\kappa}(t) \, (1-t^2)^{\lambda_\kappa-1/2}dt\not=0,\quad n=0,1,2,\ldots,
\end{equation*}
where the constant $\lambda_\kappa$ is defined in~\eqref{lambda_kappa}.
\end{teoen}

\proofen 
For $\lambda>0$, let
\begin{equation*}c_{\lambda}=\Bigr(\int\nolimits_{-1}^1 (1-t^2)^{\lambda-1/2}\,dt\Bigr)^{-1}\end{equation*}
be the normalizing constant. Then the Gegenbauer expansion of $g$ takes the form
\begin{equation}\label{Gegenbauer_expansion}
g(t)\sim\sum\limits_{n=0}^\infty b_n\frac{n+\lambda}{\lambda} C_{n}^{\lambda}(t)\quad\text{with}\quad b_n=\frac{c_{\lambda}}{C_n^{\lambda}(1)} \int\nolimits_{-1}^1 g(t) C_n^{\lambda}(t) \, (1-t^2)^{\lambda-1/2} dt,
\end{equation}
since $c_\lambda\int\nolimits_{-1}^1 \bigl(C_n^\lambda(t)\bigr)^2 \,(1-t^2)^{\lambda-1/2}\,dt=C_n^\lambda(1)\lambda/(n+\lambda)$.

For $\delta>0$, the Ces\`{a}ro $(C,\delta)$ means of the above series are
\begin{equation}\label{Cesaro_means}
S_n^\delta g(t)=\frac{1}{A_n^\delta}\sum\limits_{m=0}^n A_{n-m}^\delta b_m\frac{m+\lambda}{\lambda}C_m^\lambda(t),\quad A_n^\delta=\binom{n+\delta}n.
\end{equation}
Note that $A_{n-m}^\delta/A_n^\delta\to1$ as $n\to\infty$ for each $m$.

If $\delta>\lambda$, then it follows from \cite[Theorem~1.3]{chanillo_muckenhoupt_estimates:1993} that the sequence $\{S_n^\delta g\}$ converges uniformly on $[-1,1]$ to the function $g$.

Let $\lambda=\lambda_\kappa$ and $\delta>\lambda$. Then the sequence $\{V_\kappa(x;S_n^\delta g,y)\}$ converges uniformly in $x,y\hm\in\mathbb{S}^{d-1}$ to $V_\kappa(x;g,y)$. Indeed,
\begin{equation*}\begin{split}
\sup\limits_{x,y\in\mathbb{S}^{d-1}} \, &\bigl|V_\kappa(x;g,y)-V_\kappa(x;S_n^\delta g,y)\Bigr|\\&=\sup\limits_{x,y\in\mathbb{S}^{d-1}} \, \Bigl|\int\nolimits_{\mathbb{B}^d} g(\langle x,\xi\rangle)\, d\mu_y^\kappa(\xi)-\int\nolimits_{\mathbb{B}^d} S_n^\delta g(\langle x,\xi\rangle) \, d\mu_y^\kappa(\xi)\Bigr|\\
&\leq\sup\limits_{x,\xi\in\mathbb{S}^{d-1}} \, \bigl|g(\langle x,\xi\rangle)-S_n^\delta g(\langle x,\xi\rangle)\bigr|\\
&=\sup\limits_{t\in[-1,1]} \, \bigl|g(t)-S_n^\delta g(t)\bigr|\to0,\quad n\to\infty.
\end{split}
\end{equation*}
Using Property (II) of the $\kappa$-spherical harmonics (see Section~\ref{section_with_Dunkl_harmonic_analysis}) and \eqref{Cesaro_means}, we get
\begin{equation*}
V_\kappa(x;g,y)=\lim\limits_{n\to\infty} \, \sum\limits_{m=0}^n \frac{A_{n-m}^\delta}{A_n^\delta} \, b_m \, \sum\limits_{j=1}^{N(m,d)} S_{m,j}^\kappa(x)S_{m,j}^\kappa(y).
\end{equation*}

We can apply Theorem~\ref{main_general_theorem} with $U_n=\mathcal{A}_n^d(\kappa)$, $u_{nj}=S_{n,j}^\kappa$, and $A_{nm}=(A_{n-m}^\delta b_m)/A_n^\delta$ for $m\leq n$, $A_{nm}=0$ for $m>n$, to conclude that the condition
\begin{equation*}
\lim\limits_{n\to\infty} \, \Bigl(\frac{A_{n-m}^\delta}{A_n^\delta} \, b_m\Bigr)\not=0\quad \text{for each $m$}
\end{equation*}
is the necessary and sufficient condition for fundamentality. Obviously, this condition reduces to $b_m\not=0$ for all $m$. By \eqref{Gegenbauer_expansion}, the latter is equivalent to the integral condition described in the theorem.
\hfill$\square$

\vspace{1em}The approach used in the proof of this theorem is as in \cite[Theorem~2]{sun_cheney_article_fundamental_sets:1997}.

\section*{Acknowledgements}

This work was supported by the Russian Foundation for Basic Research (grant no.~13-01-00045) and the Ministry of Education and Science of the Russian Federation (state contract no.~1.1333.2014K).

\begin{Biblioen}

\bibitem{chanillo_muckenhoupt_estimates:1993}S. Chanillo and B. Muckenhoupt, Weak type estimates for Ces\`{a}ro sums of Jacobi polynomial series, \emph{Mem. Am. Math. Soc.} \textbf{102}:487 (1993), 1--90.

\bibitem{dai_xu_book_approximation:2013}F. Dai and Y. Xu, \emph{Approximation theory and harmonic analysis on spheres and balls}, Springer, Berlin--New York, 2013.

\bibitem{dunkl_article_reflection_groups:1988}C. F. Dunkl,  Reflection groups and orthogonal polynomials on the sphere, \emph{Math. Z.} \textbf{197} (1988), 33--60.

\bibitem{dunkl_article_operators:1989}C. F. Dunkl, Differential-difference operators associated to reflection groups, \emph{Trans. Amer. Math. Soc.} \textbf{311}:1 1989, 167--183.

\bibitem{dunkl_article_integral_kernels:1991}C. F. Dunkl, Integral kernels with reflection group invariance, \emph{Can. J. Math.} \textbf{43}:6 (1991), 1213--1227.

\bibitem{dunkl_xu_book_orthogonal_polynomials:2014}C. F. Dunkl and Y. Xu, \emph{Orthogonal polynomials of several variables}, 2nd ed., Cambridge Univ. Press, 2014.

\bibitem{Humphreys_book_reflection:1990}J. E. Humphreys, \emph{Reflection groups and Coxeter groups}, Cambridge Univ. Press, 1990.

\bibitem{rosler_article_positivity:1999}M. R\"{o}sler, Positivity of Dunkl's intertwining operator, \emph{Duke Math. J.} \textbf{98}:3 (1999), 445--463.

\bibitem{sun_cheney_article_fundamental_sets:1997}Xingping Sun and E. W. Cheney, Fundamental sets of continuous functions on spheres, \emph{Constr. Approx.} \textbf{13}:2 (1997), 245--250.

\end{Biblioen}

\noindent \textsc{Department of Applied Mathematics and Computer Science, Tula State University, Tula, Russia }

\noindent \textit{E-mail address}: \textbf{veprintsevroma@gmail.com}

\end{document}